\documentclass[a4paper,oneside]{article}
\usepackage{amsfonts}
\usepackage{amsopn}
\usepackage{amssymb}
\usepackage{euscript}
\usepackage[T2A]{fontenc}
\usepackage[ukrainian]{babel}
\usepackage[12pt]{extsizes}
\usepackage[dvips,all]{xy}
\usepackage{anysize}
\usepackage{amsmath}
\usepackage{amsthm}
\usepackage{latexsym}

\begin{document}

\begin{center}
{\large\textbf{ On existence and continuation of mild solutions of functional-differential equations of neutral type in Banach spaces }}\\[10pt]
 \end{center}

\begin{center}
Oleh Perehuda\textsuperscript{1}, Andriy Stanzhytskiy\textsuperscript{2}, Olha Martynyuk\textsuperscript{3} \\[10pt]
\small{
\textsuperscript{1} Taras Shevchenko National University of Kyiv \\  
\textsuperscript{2} Institute of Mathematics of NASU \\  
\textsuperscript{3} Yuriy Fedkovych Chernivtsi National University \\[10pt]
}
\end{center}

\underline{Abstract}. The aim of this work is to investigate the conditions for the existence and continuation of a mild solution to the initial value problem of functional-differential equations of neutral type in Banach spaces to the boundary of the domain. Based on the Schauder fixed point theorem, the existence and continuation to the boundary of the domain are proved.

\vspace{3mm}

\emph{This research was supported by NRFU project No. 2023.03/0074 Infinite-dimensional evolutionary equations
with multivalued and stochastic dynamics
}

\vspace{3mm}

1.
\underline{Introduction}.

We consider the initial value problem for the functional-differential equation of neutral type in a Banach space $X$:

\begin{equation}\label{1}
  \begin{cases}
\frac{d}{dt}(u(t)+g(t,u_{t}))=Au+f(t,u_{t}),\\
u(t)= \varphi(t), \quad t \in [-h,0].
  \end{cases}
\end{equation}

For some $h>0$ let  $C=C([-h,0];X)$ be the space of continuous functions  $\varphi : [-h,0]\rightarrow X$ with the norm  $\|\varphi\|_{C}=\sup\limits_{t\in[-h,0]} \|\varphi(t)\|_{X}$. In the sequel, the norm  $\|\varphi(t)\|_{X}$ we will denote by $\|\cdot\|_{X}$.

Let  $A$ be the infinitesimal generator of an analytic semigroup   $S(t), t\geq0$ in  $X$, $u_{t}=u(t+\theta)\in C$, for any $t\in[0,T]$.
Let  $D$ be some domain in $[-h,T]\times C$, $\partial D$ be its boundary in the sense of (\cite{Kato}, P.18), and $D=D\cup \partial D$. Mappings  $f$ and $g$ act from $D$ to $X$.

In the sequel, a solution of the problem \eqref{1} we will understand in the mild sense.

\underline{Definition 1}. A function $u:[0,T]\mapsto X$ is called a mild solution of the initial value problem  \eqref{1} on $[0,T]$ if:

1) $u(t)=\varphi (t)$, $t\in [-h,0]$;

2) $u\in C([0,T],X)$;

3) $u(t)$ satisfies the integral equation

\begin{equation}\label{2}
u(t)=S(t)(\varphi(0)+g(0,\varphi))-g(t,u_{t})-\int_{0}^{t}AS(t-s)g(s,u_{s})ds+\int_{0}^{t}S(t-s)f(s,u_{s})ds.
\end{equation}

We will study local existence of a mild solution of  \eqref{1} for  $(0,\varphi)\in D$.

A similar question in the finite-dimensional case was considered in \cite{Hal}. The author pays detailed attention to the functional-differential equation of the ordinary type. Regarding equations of the neutral type, similar results were obtained there for atomic mappings.

Regarding equations in infinite-dimensional spaces, it is worth pointing to the paper \cite{Kap} for the mild resolvability of differential equations, and the monograph \cite{Jian}, where for functional-differential equations of the ordinary type conditions for the local existence of mild solutions were obtained under the condition of continuity of the right-hand sides.

However, cases are often encountered, for example, in optimal control problems \cite{Kich}, \cite{KapN}, \cite{Lavr}, \cite{Stannos} when the conditions for the continuity of the right-hand sides are unnatural. Therefore, it is important to obtain the conditions for the existence and continuation of solutions under the conditions of the $t$-measurability  of the coefficients of the equation (such as the Carathéodory conditions). For functional-differential equations of the ordinary type, such conditions were obtained in \cite{Kich}. Note,  that even in the finite-dimensional case the solution cannot always be extended to the boundary of the domain (see the example of Myshkis, \cite{Hal}, P.60).

The aim of this work is to obtain similar results for equations of neutral type. Note that the neutral type of the equation causes the presence in the formula \eqref{2} of two new terms $g(t,u_{t})$ and $\int_{0}^{t}AS(t-s)g(s,u_{s})ds$, which significantly complicates the study of \eqref{1}.

The work consists of an introduction and three parts. The introduction describes the problem statement and a literature review. Part 2 presents the necessary concepts and formulates the main results. Part 3 is devoted to the proof of the main results. Part 4 gives an example of the obtained results for partial differential equations of parabolic type.

\vspace{3mm}
2.
\underline{Problem statement and main results}.

Let  $X$ be a reflexive Banach space, $A:X\rightarrow X$ be a linear operator with domain $D(A)$, $\sigma(-A)$ be a spectrum of $(-A)$.

\underline{Assumptions on the operator $A$}.

Assumptions (H1). \ $Re(\sigma(-A))>\delta>0$ and $A^{-1}$ is compact operator in $X$.

Then for every $\alpha\in[0,1]$ we can define fractional power  $(-A)^{\alpha}$ which is a closed linear operator with domain $D((-A)^{\alpha}$ (see \cite{Pazy} for details). let us denote by $X_{\alpha}$ the Banach space $D((-A)^{\alpha})$ supplied with the norm  $\|u\|_{\alpha}:=\|(-A)^{\alpha}u\|$, which is equivalent to the graph norm of $(-A)^{\alpha}$. Let us denote $X_{0}=X$. From \cite{Henry}, Sec.1.4 we deduce that if $A^{-1}$ is a compact operator, then the semigroup  $S(t)$ is compact for $t>0$.  Then under assumption (H1)  Theorem 3.2 of \cite{Pazy} implies that the semigroup  $S(t)$ is continuous in the uniform operator topology for $t>0$.
Thus, due to (\cite{Pazy}, Th.3.3), we can conclude that the operator  $A$  has compact resolvent. The last property guarantees the following result (\cite{Henry},Th.1.4.8):

\textbf{Lemma 1.}  \emph{Under assumption  (H1) the embedding  $X_{\alpha}\subset X_{\beta}$ is compact for  $0\leq \beta<\alpha\leq1$}.

\textbf{Lemma 2.}  \emph{(\cite{Henry},Th.1.4.3). Under assumption  (H1) for every $\alpha\geq 0$ there exists  $C_{\alpha}>0$ such that
$$\|(-A)^{\alpha}S(t)\|\leq C_{\alpha}t^{-\alpha}\exp^{-\delta t}, \ \ t>0.$$}

In particular,  $\|S(t)\|\leq C_{0}\exp^{-\delta t}$, \ $t>0$.

\underline{Assumptions on nonlinear terms}.

Let  $U$ be an open subset of  $C$ and  $f:[0,T]\times U\rightarrow X$.

Assumption (H2).

1) for every $t\in[0,T]$ the mapping $f$ is continuous with respect to $\varphi$;

2) for every  $\varphi \in U$ the mapping $f$ is measurable with respect to $t$;

3) for any $R>0$ there exists integrable function $m_{R}(t)\in L^{p}(0,T)$, $p>1$ such that  $\|f(t,\varphi)\|\leq m_{R}(t)$, if $\|\varphi\|\leq R$.

4) there exist  $\alpha\in [0,1]$ and $M_{g}\in(0,1)$ such that  $g:[0,T]\times U\rightarrow X_{\alpha}$;

5) $\|g(t,\varphi_{1})-g(t,\varphi_{2})\|_{\alpha}\leq M_{g}\|\varphi_{1}-\varphi_{2}\|_{C}$;

6)  function $g$ is continuous with respect to $t$ in $X_{\alpha}$ uniformly over $\varphi\in U$.

\textbf{Theorem  1.}(local existence)   \emph{Suppose that assumptions (H1)-(H3) hold. Then for every  $\varphi\in U$ there exist  $t_{1}=t_{1}(\varphi)\in [0,T]$ and continuous function  $u\in [-h,t_{1}]\rightarrow X$ such that  $u(t)$ is a mild solution of  \eqref{1} on  $[0,t_{1}]$ with the initial function  $\varphi$}.

Note that this theorem takes place for every initial point  $t_{0}\in (0,T)$. So, we can continue the solution on the maximal interval  $[t_{0}-h,t_{0}]$.

We consider  $D=(0,T)\times U$, its boundary  $\partial D$, and its closure  $\overline{D}=D\bigcup \partial D$.

\textbf{Theorem 2.}(continuation)   \emph{Assume that conditions of Theorem 1 are valid. Then the solution with initial data
 $(t_{0},\varphi)\in D$ (i.e., on the interval $[t_{0}-h,t_{0}]$) exists on the maximal interval  $[t_{0},\tau)$, $\tau >t_{0}$, and  $(\tau, u_{\tau})\in \partial D$}.
%

Further, we will use the following fixed point theorem.

\textbf{Theorem} (Krasnoselskii \cite{Bur})  \emph{Let  $M$ be a non-empty closed convex subset of a Banach space $(S,\|\cdot\|)$. Assume that operators  $A$ and $B$ act from $M$ into  $S$ such that:}

\emph{(i) $Ax+By\in M \ \forall \ x,y\in M$};

\emph{(ii) $A$ is a continuous operator and  $AM$ belongs to a compact set};

\emph{(iii) $B$ is a contraction operator with the constant  $\alpha<1$}.

\emph{Then there exists  $y\in M$ such that  $Ay+By=y$.}

3.
\underline{Proofs of the main results}.

The proof of Theorem 1.

We want to apply the Krasnoselskii fixed point theorem. It should be noted that due to assumptions on $A$  there exists  $M>0$ such that  $$\|S(t)\|\leq M \ \forall \ t\in [0,T].$$ We choose  $\beta>0$ such that
$$B_{\beta}=\{\psi \in C: \|\varphi-\psi\|_{C}\leq \beta\}\subset U.$$
Analogously to  \cite{Hal}, \cite{Jian} we introduce closed convex bounded set
$$\mathcal{A}(\delta, \beta)=\{y\in C([-h,\delta];X):y_{0}=\varphi, y_{t}\in B_{\beta}, t\in [0,\delta]\},$$ where $\delta$ is sufficiently small.  On this set we consider the operator
\begin{equation}\label{3}
G(y)(t)=
\begin{cases}\varphi(t),  t\in [-h,0],\\
S(t)(\varphi(0)+g(0,\varphi))-g(t,y_{t})-\int_{0}^{t}AS(t-s)g(s,y_{s})ds+\int_{0}^{t}S(t-s)f(s,y_{s})ds, t\in [0,\delta].
  \end{cases}
  \end{equation}
We split this mapping as  $G(y)=\Phi (y)+\Psi (y)$, where
\begin{equation}\label{4}
\Phi(y)(t)=
\begin{cases}\varphi(t),  t\in [-h,0],\\
S(t)(\varphi(0)+g(0,\varphi))-\int_{0}^{t}AS(t-s)g(s,y_{s})ds+\int_{0}^{t}S(t-s)f(s,y_{s})ds, t\in [0,\delta].
  \end{cases}
\end{equation}
and
$$\Psi(y)(t)=
\begin{cases}0,  t\in [-h,0],\\
-g(t,y_{t}), t\in [0,\delta].
  \end{cases}
 $$
Let us check for the given mappings whether the conditions of the Krasnoselskii theorem are fulfilled.

Let us show the continuity of  $\Phi(y)$.

Let $y^{(n)}(t)\in \mathcal{A}(\delta,\beta)$ and  $\sup\limits_{t\in [0,\delta]}\|y^{(n)}(t)-y(t)\|\rightarrow 0, n\rightarrow\infty$.

Then
$$
\Phi(y^{(n)})(t)=
\begin{cases}\varphi(t),  t\in [-h,0],\\
S(t)(\varphi(0)+g(0,\varphi))-\int_{0}^{t}AS(t-s)g(s,y_{s}^{(n)})ds+\int_{0}^{t}S(t-s)f(s,y_{s}^{(n)})ds,
  \end{cases}
$$
$$
t\in [0,\delta].
$$
So,
\begin{equation}\label{5}
\begin{array}{l}
\sup\limits_{t\in [-h,\delta]}\|\Phi(y^{(n)}(t)-\Phi(y)(t)\|\leq \sup\limits_{t\in [0,\delta]}\|\int_{0}^{t}AS(t-s)(g(s,y_{s}^{(n)}-g(s,y_{s}))ds\|+ \\
\linebreak
+\sup\limits_{t\in [0,\delta]}\|\int_{0}^{t}S(t-s)(f(s,y_{s}^{(n)}-f(s,y_{s}))ds\|.
\end{array}
\end{equation}
For the first summand we have:
\begin{equation}\label{6}
\begin{array}{l}
\|\int_{0}^{t}AS(t-s)(g(s,y_{s}^{(n)}-g(s,y_{s}))ds\|\leq \|\int_{0}^{t}(-A)^{1-\alpha}S(t-s)\|\|(-A)^{\alpha}(g(s,y_{s}^{(n)}-g(s,y_{s}))\|ds\leq \\
\\
\leq \int_{0}^{t}\|(-A)^{1-\alpha}S(t-s)\|ds \sup\limits_{s\in [0,\delta]}\|g(s,y_{s}^{(n)}-g(s,y_{s}\|_{X_{\alpha}}\leq \\
\\
\leq\int_{0}^{t}C_{1-\alpha }t^{\alpha-1}ds M_{g} \sup\limits_{s\in [0,\delta]}\|y_{s}^{(n)}-y_{s}\|=C_{1-\alpha}\sup\limits_{s\in [0,\delta]}\|y_{s}^{(n)}-y_{s}\|\rightarrow 0, n\rightarrow\infty.
\end{array}
\end{equation}
It should be noted that  $\sup\limits_{t\in [0,\delta]}\|y_{t}^{(n)}-y_{t}\|_{C}\leq \sup\limits_{t\in [0,\delta]}\|y_{t}^{(n)}-y_{t}\|, n\rightarrow\infty.$

Then for every  $t\in[0,\delta]$ we get  $f(t,y_{t}^{(n)}\rightarrow f(t,y_{t}, \ n\rightarrow\infty$ due to the continuity of  $\varphi\mapsto f(t,\varphi)$.

Then  the definition of  $\mathcal{A}(\delta,\beta)$ and Lebesgue's dominated convergence theorem yield
\begin{equation}\label{7}
\sup\limits_{t\in [0,\delta]}\|\int_{0}^{t}S(t-s)(f(s,y_{s}^{(n)}-f(s,y_{s}))ds\|\leq \int_{0}^{\delta}M\|f(s,y_{s}^{(n)}-f(s,y_{s}\|ds\rightarrow 0.
\end{equation}

Let us prove that  $\Phi(\mathcal{A}(\delta,\beta))$ is contained in a compact set. To this end, we will use the following statements.

\textbf{Proposition 1.}(\cite{DZ}, Prop.8.4) \emph{ If  $S(t), t>0$ is a compact operator then for any $p\geq 1$ the mapping   $$F(f)(t)=\int_{0}^{t}S(t-s)f(s)ds, \ f\in L^{p}(0,T;X), \ t\in [0,T]$$  is a compact operator from  $L^{p}(0,T;X)$ to  $C([0,T];X)$.}

\textbf{Proposition 2.} \cite{Stan1} \emph{ Under assumptions on the operator $A$ the mapping $$(Bz)(t)=\int_{0}^{t}AS(t-s)z(s)ds$$ is a compact operator from  $C([0,T],X_{\alpha})$ to $C([0,T],X)$}.

Further, we will use the infinite dimensional analogous of the Arzela-Askoli theorem. We need to prove that

1) for every  $t\in [0,\delta]$ the set  $\{\Phi(y)(t):y\in\mathcal{A}(\delta,\beta)\}$ is precompact in  $X$;

2) the family of functions $\Phi(y)(t)$ is equicontinuous, i.e., for any $\epsilon >0$ there exists  $\sigma >0$ such that  $\|\Phi(y)(t_{1})-\Phi(y)(t_{2})\|\leq \varepsilon$, if  $y\in \mathcal{A}(\delta,\beta), |t_{1}-t_{2}|<\sigma, t_{1}, t_{2}\in [0,\delta]$.

Note that from (H3) we get for $t\in[0,\delta]$
\begin{equation}\label{8}
\sup\limits_{t\in [0,\delta]}\|g(t,y_{t})\|_{\alpha}\leq \sup\limits_{t\in [0,\delta]}\|(g(t,y_{t}^{(n)})-g(t,\varphi))\|_{\alpha}+\|g(t,\varphi)\|_{\alpha}\leq M_{g}\beta+\|g(t,\varphi)\|_{\alpha}\leq C_{1}.
\end{equation}

Then the compactness of the first integral summand in  \eqref{4} is the consequence of Proposition 2.

Further, due to the definition of  $\mathcal{A}(\delta,\beta)$, we get
 $$\sup\limits_{t\in [0,\delta]}\|y(t)\|_{C}\leq \sup\limits_{t\in [0,\delta]}\|y_{t}-\varphi\|_{C}+\|\varphi\|_{C}\leq\beta+\|\varphi\|_{C}.$$

Using condition 3) from (H2), we obtain that  $\|f(t,y_{t})\|\leq m_{R}(t)$ for $R=\beta+\|\varphi\|_{C}$. Then the set  $f(t,y_{t})$
is bounded in  $L^{p}(0,T;X)$. Therefore, Proposition 1 yields the compactness of the second integral summand in  \eqref{4}. Thus, property 1) is proved.
Now let us prove  the equicontinuity property. let us fix an arbitrary  $\varepsilon$. We note that  the continuity of  $S(t)(\varphi(0)+g(0,\varphi))$,  $t\in [0,\delta]$, implies its uniform continuity on $[0,\delta]$. Therefore, there exists  $\sigma_{1}$ such that for all $t_{1},t_{2}\in [0,\delta]$ with $|t_{1}-t_{2}|<\sigma_{1}$ we have
\begin{equation}\label{9}
\|(S(t_{1})-S(t_{2}))(\varphi(0)+g(0,g(0,\varphi))\|<\frac{\varepsilon}{3}.
\end{equation}

After that, using Lemma 2 and estimate  \eqref{8}, we get
\begin{equation}\label{10}
\begin{array}{c}
\|\int_{t_{1}}^{t_{2}}AS(t-s)g(s,y_{s})ds\|\leq \int_{t_{1}}^{t_{2}}\|(-A)^{1-\alpha}S(t-s)\|ds\sup\limits_{t\in [0,\delta]}\|g(s,y_{s}\|\leq \\
\\
\leq \frac{C_{1-\alpha}}{\alpha}|t_{2}-t_{1}|C_{1}<\frac{\varepsilon}{3},
\end{array}
\end{equation}
if  $|t_{1}-t_{2}|<\sigma_{2}$.

Due to the absolutely continuity of Lebesgue's integral we deduce that
 \begin{equation}\label{11}
 \|\int_{t_{1}}^{t_{2}}S(t-s)f(s,y_{s})ds\|\leq M\int_{t_{1}}^{t_{2}}m_{R}(t)<\frac{\varepsilon}{3},
 \end{equation}
 if $|t_{2}-t_{1}|<\sigma_{3}$, and $R=\beta+\|\varphi\|_{C}$.

 Then, choosing  $\sigma=min\{\sigma_{1},\sigma_{2},\sigma_{3}\}$, from  \eqref{9}-\eqref{11} we deduce equicontinuity of the family  $\Phi(y)(t)$. As  $\Phi(\mathcal{A}(\delta,\beta))$ is a compact set in  $C([-h,0],X)$, so we have that assumption (ii) of the Krasnoselskii theorem is fulfilled. Let us verify assumption (iii).

For every $y^{(1)}(t)$ and  $y^{(2)}(t)$, which belongs to $\mathcal{A}(\delta,\beta)$, we have
$$
\begin{array}{l}
\sup\limits_{t\in [-h,\delta]}\|\Phi(y^{(1)}(t))-\Phi(y^{(2)}(t))\| = \sup\limits_{t\in [0,\delta]}\|g(t,y_{t}^{(1)})-g(t,y_{t}^{(2)})\|\leq \\
\leq\sup\limits_{t\in [0,\delta]}\|g(t,y_{t}^{(1)})-g(t,y_{t}^{(2)})\|_{\alpha}
\leq M_{g}\sup\limits_{t\in [0,\delta]}\|y_{t}^{(1)}-y_{t}^{(2)}\|=M_{g}\sup\limits_{t\in [-h,\delta]}\|y_{t}^{(1)}-y_{t}^{(2)}\|
\end{array}
$$
Because of the inequality  $M_{g}<1$ we get (iii).

Let us prove (i). For arbitrary $y\in \mathcal{A}(\delta,\beta)$ and  $z\in \mathcal{A}(\delta,\beta)$ we get
\begin{equation}\label{12}
\Phi(y)(t)+\Psi(z)(t)=
\begin{cases}\varphi(t),  t\in [-h,0],\\
S(t)(\varphi(0)+g(0,\varphi))-g(t,z_{t})-\int_{0}^{t}AS(t-s)g(s,y_{s})ds+\int_{0}^{t}S(t-s)f(s,y_{s})ds.
  \end{cases}
\end{equation}
Then for every  $t\in [0,\delta]$ we have
\begin{equation}\label{13}
\begin{array}{c}
\|\Phi(y)_{t}+\Psi(z)_{t}-\varphi_{0}\|=\sup\limits_{\theta\in [-h,0]}\|\Phi(y)_{t}+\Psi(y)_{t}-\varphi(\theta)\|\leq \\
\leq\sup\limits_{\theta\in [-h,-t]}\|\Phi(y)(t+\theta)+\Psi(y)(t+\theta)-\varphi(\theta)\|+\sup\limits_{\theta\in[-t,0]}\|\Phi(y)(t+\theta)+\Psi(y)(t+\theta)-\varphi(\theta)\|.
\end{array}
\end{equation}
From the uniform continuity of $\varphi$ we deduce
\begin{equation}\label{14}
\sup\limits_{\theta\in[-h,-t]}\|\Phi(y)(t+\theta)+\Psi(y)(t+\theta)-\varphi(\theta)\|= \sup\limits_{\theta\in[-h,-t]}\|\varphi(t+\theta)-\varphi(\theta)\|\leq\frac{\beta}{k}
\end{equation}
for all $t\in[0,\delta]$ and for sufficiently small $\delta$. Here $k$ is some natural number, which we choose latter, and  $\delta=\delta(k)$. Further, for  $t+\theta\in[0,\alpha]$ we get
\begin{equation}\label{15}
\begin{array}{l}
\Phi(y)(t+\theta)+\Psi(z)(t+\theta)-\varphi(\theta)=S(t+\theta)(\varphi(0)+g(0,\varphi))-g(t+\theta,z_{t+\theta})-\\
\\
-\int_{0}^{t+\theta}AS(t+\theta-s)g(s,y_{s})ds +\int_{0}^{t+\theta}AS(t+\theta-s)f(s,y_{s})ds-\varphi(\theta)= \\
\\
=S(t+\theta)\varphi(0)-\varphi(0)-\varphi(\theta)+\varphi(0)+S(t+\theta)g(0,\varphi)-g(t+\theta,z_{t+\theta})-\\
\\
-\int_{0}^{t+\theta}AS(t+\theta-s)g(s,y_{s})ds +\int_{0}^{t+\theta}AS(t+\theta-s)f(s,y_{s})ds.
\end{array}
\end{equation}
Here  $\theta\in[-\alpha,0]$. Then
\begin{equation}\label{16}
\|S(t+\theta)\varphi(0)-\varphi(0)\|\leq\frac{\beta}{k}, t\in[0,\delta], \delta=\delta(k),
\end{equation}
due to the  $C_{0}$-continuity of  $S(t)$.

Choosing  $\delta$ (depending on $k$) sufficiently small, we get
\begin{equation}\label{17}
\|\varphi(\theta)-\varphi(0)\|\leq\frac{\beta}{k}.
\end{equation}
Further,
\begin{equation}\label{18}
S(t+\theta)g(0,\varphi)-g(t+\theta,z_{t+\theta})=S(t+\theta)g(0,\varphi)-g(0,\varphi)+g(0,\varphi)-g(t+\theta,z_{t+\theta}).
\end{equation}
But
\begin{equation}\label{19}
\|S(t+\theta)g(0,\varphi)-g(0,\varphi)\|\leq\frac{\beta}{k}, t\in[0,\delta(k)].
\end{equation}
We have
\begin{equation}\label{20}
\begin{array}{l}
\|g(t+\theta,z_{t+\theta})-g(0,\varphi)\|\leq\|g(t+\theta,z_{t+\theta})-g(t+\theta,\varphi)\|_{\alpha}+\|g(t+\theta,\varphi)-g(0,\varphi)\|_{\alpha}\leq \\
\\
\leq M_{g}\|z_{t+\theta}-\varphi\|_{C}+\|g(t+\theta,\varphi)-g(0,\varphi)\|_{\alpha}\leq M_{g}\beta+\frac{\beta}{k}.
\end{array}
\end{equation}
The last estimate follows from the definition of  $\mathcal{A}(\delta,\beta)$ and continuity of the mapping  $g$ with respect to $t$.

Let us estimate the first integral summand in \eqref{15}.

We have
\begin{equation}\label{21}
\begin{array}{l}
\|\int_{0}^{t+\theta}AS(t+\theta-s)g(s,y_{s})ds\|\leq \int_{0}^{t+\theta}C_{1-\alpha}(t+\theta-s)^{\alpha-1}ds\sup\limits_{s\in[0,t+\theta]}(M_{g}\|y_{s}-\varphi\|+\|g(s,\varphi)\|_{\alpha})\leq
\\
\leq\frac{C_{1-\alpha}}{\alpha}(t+\theta)^{\alpha}(M_{g}\beta+\sup\limits_{s\in[0,\alpha]}\|g(s,\varphi)\|_{\alpha})\leq\frac{\beta}{k}
\end{array}
\end{equation}
for sufficiently small  $\delta=\delta(k)$.

To estimate the last integral in  \eqref{15} we have
\begin{equation}\label{22}
\|\int_{0}^{t+\theta}S(t+\theta-s)f(s,y_{s})ds\|\leq M\int_{0}^{\delta}m_{\beta}(s)\leq\frac{\beta}{k},
\end{equation}
due to the absolutely continuity of Lebesgue's integral.

Thus, from  \eqref{13}, \eqref{14}, \eqref{17}-\eqref{22} we deduce that
\begin{equation}\label{23}
\sup\limits_{t\in[0,\delta(k)]}\|\Phi(y)_{t}+\Psi(z)_{t}-\varphi_{0}\|\leq\frac{7\beta}{k}+M_{g}\beta.
\end{equation}

Now let us choose  $k$ such that  $\frac{7}{k}+M_{g}<1$. Further, we choose  $\delta$ such that inequalities  \eqref{14}, \eqref{17}-\eqref{22} are fulfilled. Taking into account continuity of  $\Phi$ and $g(t,\varphi)$, we get from \eqref{23} condition (i) of the Krasnoselskii theorem. After applying this theorem, we deduce that there exists  $y\in\mathcal{A}(\delta,\beta)$ such that  $G(y)=y$. Theorem 1 is proved.

The proof of Theorem 2.

We need the following lemmas.

\textbf{Lemma 3.1.} \emph{ Under conditions of Theorem 1 if  $W\subset D$ and  $W$ is a compact set, then there exists  $\delta>0$ such that for arbitrary initial data  $(t_{0},\varphi)\in W$ the solution of  \eqref{1} exists on the interval  $[t_{0},T_{0}+\delta]$.}

\underline{The proof of Lemma 3.1.} As  $W$ is a compact set, so there exists an open set  $V$ such that  $W\subset V\subset D$. Therefore, there exists a function  $m(t)$ such that the inequality 3) takes place simultaneously for all  $(t_{0},\varphi)\in W$.

Compactness of $W$ implies the inequality  \eqref{8} for all $\varphi$ with $(t_{0},\varphi)\in W$. By standard considerations related to the existence for every  $\varepsilon >0$  a finite $\varepsilon$-net in  $W$, we can get that the inequalities  \eqref{14}, \eqref{17}-\eqref{22} take place uniformly over all  $\varphi$ with  $(t_{0},\varphi)\in W$. So, there exists  $\delta >0$ such that all mentioned above inequalities take place on  $[0,\delta]$ uniformly over  $(t_{0},\varphi)\in W$. The rest of the proof is the consequence of Theorem 1.

\textbf{Lemma 3.2.} \emph{If  $u(t)$ is a non-continuing solution of  \eqref{1} on $[0,\tau]$, then for every compact set  $W\subset D$ there exists $t_{W}$ such that  $(t,u_{t})\not \in W$ for $t\in[t_{W},\tau)$.}

\underline{The proof of Lemma 3.2.} Since $W$ is a compact set in $D$, so Lemma 1 implies that for every  point $(c,\varphi)$ equation  \eqref{1} has a solution  $u(t)$ such that  $u(t)=\varphi(t)$ for $t\in [c-h,c]$, which exists at least on the interval  $[c,c+\delta]$. Arguing by contradiction, we can find a sequence  $t_{k}\rightarrow \tau$, $k\rightarrow\infty$ and $\psi\in C$ such that  $(t_{k},u_{t_{k}})\in W$, $(\tau,\psi)\in W$, $(t_{k},u_{t_{k}})\rightarrow (\tau,\psi)$ as $k\rightarrow\infty$. The last statement is a consequence of the  compactness of $W$.

So, for arbitrary  $\varepsilon\in(0,h)$ we have  $$\lim\limits_{k\rightarrow\infty}\sup\limits_{\theta\in [-h,-\varepsilon]}\|u(t_{k}+\theta)-\psi(\theta)\|=0.$$ Therefore,  $$u(t_{k}-\varepsilon)\rightarrow \psi(-\varepsilon), k\rightarrow\infty.$$ Thus,  $u(\tau-\varepsilon)=\psi(-\varepsilon)$ and we can conclude that there exists $\lim\limits_{t\rightarrow\tau-0}u(t)$ and  $u(t)$ can be continued  to continuous function on  $[-h,\tau]$, if we put  $u(\tau)=\psi(0)$. But  $(\tau,u_{\tau})\in W$. So, due to lemma 1 there exists a solution with the initial data  $(\tau,u_{\tau})$ on the interval  $[\tau,\tau+\delta]$, which is a contradiction. Lemma is proved.

Now we are in position to prove Theorem 2.

We have that the solution of  \eqref{1} with initial data $(t_{0},\varphi)\in D$ exists on some interval  $[t_{0},t_{0}+\delta]$. Due to Zorn's lemma this solution can be continued on the maximal interval  $(t_{0},\tau)$. Let  $K$ be a closed bounded set in $D$. Let us show that there exists  $t_{k}$ such that  $(t,u_{t})\not\in K$ for  $t\in (t_{k},\tau)$. Assume the contrary. Then there exists a sequence $t_{k}\rightarrow\tau-0$ such that  $(t_{k},u_{t_{k}})\in K$  for all $k$. Then for every  $\varepsilon>0$ the continuity of the solution  $u(t)$ in the space  $X$  on $[t_{0}-h,\tau-\varepsilon]$ implies that uniform continuity of  $u_{t}$ in the space $C$ for  $t\in[t_{0},\tau-\varepsilon]$. Then the set $P=\{(t,u_{t}):t\in[t_{0},\tau)\}$ is bounded. Indeed, if it is not bounded, then there exists a sequence  $(s_{k},u_{s_{k}})$  such that  $s_{k}\rightarrow\tau-0$, and  $\|u_{s_{k}}\|_{C}\rightarrow\infty, k\rightarrow\infty$. But $(t_{k},u_{t_{k}})\in K$, so  $u_{t_{k}}$ are bounded in $C$, which contradicts the uniform continuity of $u_{t}$ on $[t_{0},\tau-\varepsilon]$. Therefore,  $P$ is bounded and  $\overline{P}$ belongs to $K$.

Let us show that  $P$ belongs to some compact set in  $D$. It is sufficient to prove precompactness of the set  $R=\{u_{t}, t\in[t_{0},\tau)\}$. According to the Arzela-Askoli theorem we need to prove that

1) for every  $\theta\in[-h,0]$ the set  $R(\theta)=\{u(t+\theta), t\in[t_{0},\tau)\}$ is precompact in $X$;

2) the family of functions $\{u(t+\theta), \theta\in[-h,0], t\in[t_{0},\tau)\}$ is equicontinuous.

For proving the precompactness of $R(\theta)$ we will construct a finite $\delta$-net for every $\delta>0$.

Without loss of generality we assume that  $t_{0}=\theta$. Let us fix sufficiently small $\mu<\tau$. We split the set $R(\theta)$ into two parts $R(\theta)=R_{1}(\theta)\cup R_{1}(\theta)$:
$$
R_{1}(\theta)=\{u(t+\theta), t+\theta\in[-h,\mu]\}, \  R_{2}(\theta)=\{u(t+\theta), t+\theta\in(\mu,\tau)\}.
$$
Uniform continuity of $u(t+\theta)$ on $[-h,\mu]$ implies the existence of a finite $\delta$-net for  $R_{1}(\theta)$.

Let us consider $R_{2}(\theta)$. As  $t+\theta>0$, then for elements of $R_{2}(\theta)$ we have
\begin{equation}\label{24}
\begin{array}{l}
u(t+\theta)=S(t+\theta)(\varphi(0)+g(0,\varphi))-g(t+\theta,u_{t+\theta})-\int_{0}^{t+\theta}AS(t+\theta-s)g(s,u_{s})ds+
\\
\\
+\int_{0}^{t+\theta}S(t+\theta-s)f(s,u_{s})ds.
\end{array}
\end{equation}
For every  $\varepsilon\in(0,\mu)$ we consider the set $R_{\varepsilon}(\theta)$, whose elements have the following representation:
\begin{equation}\label{25}
\begin{array}{l}
u^{(\varepsilon)}(t+\theta)=S(t+\theta)(\varphi(0)+g(0,\varphi))-g(t+\theta,u_{t+\theta})-\int_{0}^{t+\theta-\varepsilon}AS(t+\theta-s)g(s,u_{s})ds+
\\
\\
+\int_{0}^{t+\theta-\varepsilon}S(t+\theta-s)f(s,u_{s})ds.
\end{array}
\end{equation}
As  $\overline{P}\subset K$ is bounded, so due to assumption 3) there exists an integrable function  $m(t)$ such that
\begin{equation}\label{26}
|f(t,u_{t})|\leq m(t), t\in[0,\tau)
\end{equation}
and
\begin{equation}\label{27}
\begin{array}{l}
\sup\limits_{s\in [0,\tau)}\|g(s,u_{s})\|_{\alpha}\leq \sup\limits_{s\in [0,\tau)}\|g(s,u_{s}-g(s,\varphi)\|_{\alpha}+\sup\limits_{s\in [0,T]}\|g(s,\varphi)\|_{\alpha}\leq
\\
\\
\leq M_{g}(\sup\limits_{s\in [0,\tau)}\|u_{s}\|+\|\varphi\|)+\sup\limits_{s\in [0,T]}\|g(s,\varphi)\|\leq C,
\end{array}
\end{equation}
for some $C>0$.

Applying the semigroup property, we get
\begin{equation}\label{28}
\begin{array}{l}
u^{(\varepsilon)}(t+\theta)=S(\varepsilon)(\varphi(0)+g(0,\varphi))-g(t+\theta,u_{t+\theta})-S(\varepsilon)\int_{0}^{t+\theta-\varepsilon}AS(t+\theta-\varepsilon-s)g(s,u_{s})ds+
\\
\\
+S(\varepsilon)\int_{0}^{t+\theta-\varepsilon}S(t+\theta-s-\varepsilon)f(s,u_{s})ds.
\end{array}
\end{equation}
Since  $t+\theta-\varepsilon>0$, then
$$
\|S(t+\theta-\varepsilon)(\varphi(0)+g(0,\varphi))\|\leq M\|\varphi(0)+g(0,\varphi)\|.
$$
As the set  $S(t+\theta-\varepsilon)(\varphi(0)+g(0,\varphi))$ is uniformly bounded for  $t\in[0,\tau)$, then the set $S(\varepsilon)S(t+\theta-\varepsilon)(\varphi(0)+g(0,\varphi))$ is precompact in  $X$.

Let us prove precompactness of the set $\{g(t+\theta, u_{t+\theta}), t+\theta\in[\mu,\tau)\}$ in $X$. From \eqref{27} and compact embedding  $X_{\alpha}\subset X$ we derive that the set  $G(t_{1})=\{g(t_{1}, u_{t+\theta}), t+\theta\in[\mu,\tau)\}$ is precompact in $X$ for every  $t_{1}\in[0,T]$. So, for every $t_{1}\in[0,T]$ it has a finite $\varepsilon$-net  $\{z_{1},...,z_{p}\}, p=p(\varepsilon,t_{1})$.

Let us construct such a net for the set $G=\{g(s, u_{t+\theta}), t+\theta\in[\mu,\tau), s\in[0,T]\}$.

Due to the continuity with respect to $t\in[0,T]$ of the mapping $g(t,\varphi)$, uniformly over  $\varphi\in U$, we deduce that for every  $\varepsilon>0$ there exist a finite set  $\{t_{1}(\varepsilon), ..., t_{r}(\varepsilon)\}, r=r(\varepsilon)$ and  $\delta>0$ such that if  $|t-t_{i}|<\delta$ then
\begin{equation}\label{29}
\|g(t,\varphi)-g(t_{i},\varphi)\|<\varepsilon.
\end{equation}

Since each set  $G(t_{i}(\varepsilon))$ is precompact, so it has a finite  $\varepsilon$-net. Therefore, the union ovet  $t_{i}$ of these $\varepsilon$-nets will be $2\varepsilon$-net for the set  $G$. Then the set  $G$  is precompact in $X$ and we can conclude that the set $\{g(t+\theta, u_{t+\theta}), t+\theta\in[\mu,\tau)\}$ is precompact in  $X$.

From \eqref{25} we get
\begin{equation}\label{30}
\begin{array}{l}
\|\int_{0}^{t+\theta-\varepsilon}S(t+\theta-\varepsilon-s)f(s,u_{s}ds\|\leq \int_{0}^{t+\theta-\varepsilon}\|S(t+\theta-\varepsilon-s)\|\|f(s,u_{s}\|ds\leq
\\
\\
\leq M\int_{0}^{t+\theta-s}m(s)ds\leq \int_{0}^{T}m(s)ds\leq C_{1},
\end{array}
\end{equation}
for some positive constant  $C_{1}$.

Then, since  $S(\varepsilon)$ is a compact operator, we conclude that the set
$$\{S(\varepsilon)\int_{0}^{t+\theta-s}S(t+\theta-\varepsilon-s)f(s,u_{s})ds,   t+\theta\in[\mu,\tau)\}$$
 is precompact in  $X$.

In the same way, using \eqref{21} and  \eqref{27}, we  show that the set
$$\{S(\varepsilon)\int_{0}^{t+\theta-s}AS(t+\theta-\varepsilon-s)g(s,u_{s})ds,   t+\theta\in[\mu,\tau)\}$$
is precomapct in $X$.

So, for every $\varepsilon\in(0,\mu)$ the set  $R_{\varepsilon}(\theta)$ is precompact in $X$, therefore, for every  $\delta>0$ it has a finite  $\delta$-net. Let us denote it by  $\{u^{(\varepsilon)}(t_{1}+\theta),...,u^{(\varepsilon)}(t_{p}+\theta)\}$, $p=p(\delta)$.

From \eqref{24} and \eqref{25} we get
\begin{equation}\label{31}
\|u(t+\theta)-u^{(\varepsilon)}(t+\theta)\|\leq \|\int_{t+\theta-\varepsilon}^{t+\theta}AS(t+\theta-s)g(s,u_{s})ds\|+\|\int_{t+\theta-\varepsilon}^{t+\theta}S(t+\theta-s)f(s,u_{s})ds\|.
\end{equation}
Moreover,
\begin{equation}\label{32}
\begin{array}{l}
\|\int_{t+\theta-\varepsilon}^{t+\theta}AS(t+\theta-s)g(s,u_{s})ds\|\leq \int_{t+\theta-\varepsilon}^{t+\theta}C_{1-\alpha}(t+\theta-s)^{\alpha-1}ds\sup\limits_{s\in[0,\tau)}\|g(s,u_{s}\|=
\\
\\
=\frac{C_{1-\alpha}}{\alpha}\varepsilon^{\alpha}C_{1}\rightarrow 0, \varepsilon\rightarrow 0,
\end{array}
\end{equation}
uniformly over  $t+\theta$.

Analogously,
\begin{equation}\label{33}
\|\int_{t+\theta-\varepsilon}^{t+\theta}S(t+\theta-s)f(s,u_{s})ds\|\leq \int_{t+\theta-\varepsilon}^{t+\theta}m(s)ds\rightarrow 0,
 \varepsilon\rightarrow 0,
\end{equation}
uniformly over  $t+\theta$. Finally,
\begin{equation}\label{34}
\sup\limits_{t+\theta\in[\mu,\tau)}\|u^{\varepsilon}(t+\theta)-u(t+\theta)\|\rightarrow 0, \varepsilon\rightarrow 0.
\end{equation}

Using the  $\delta$-net for $R_{\varepsilon}(\theta)$, we can construct $3\sigma$-net for  $R_{2}(\theta)$:
\begin{equation}\label{35}
\{u(t_{1}+\theta),...u(t_{p}+\theta)\},
\end{equation}
where we choose elements $u(t_{i}+\theta)$ from $R_{2}(\theta)$ such that  $\sup\limits_{t\in(0,\tau)}\|u(t+\theta)-u_{\varepsilon}(t+\theta)\|<\delta$ (it is possible due to  \eqref{33}). Therefore, the finite $\delta$-net both for  $R_{1}(\theta)$ and $R_{2}(\theta)$  has been constructed. It means that  $R(\theta)$ is precomapct in $X$, and property 1) is proved.

Let us show equicontinuity of the family  $\{u(t+\theta), \theta\in[-h,0], t\in[0,\tau)\}$. We need to prove that for any $\varepsilon>0$ there exists  $\delta>0$ such that for all $t\in[0,\tau)$
\begin{equation}\label{36}
\|u(t+\theta_{1})-u(t+\theta_{2})\| < \varepsilon
\end{equation}
for all $\theta_{1},\theta_{2}\in[-h,0]$ with $|\theta_{1}-\theta_{2}|<\delta$.

Let us put  $\theta_{2}=\theta_{1}+r$. For sufficiently small  $\mu\in(0,\tau)$ we consider three cases.

1. $t+\theta_{2}\leq\mu$. In this case equicontinuity of $u(t+\theta)$ follows from the uniform continuity of  $u(t)$ on $[-h,\mu]$.

2. $t+\theta_{1}<\mu, t+\theta_{2}\geq \mu$. In this case we have
$$
\|u(t+\theta_{2})-u(t+\theta_{1})\|\leq \|u(t+\theta_{2})-u(\mu)\|+\|u(t+\theta_{1})-u(\mu)\|.
$$
Choosing sufficiently small $\delta$ we can guarantee that the points $t+\theta_{1}$ and $t+\theta_{2}$ belongs to small left and right neighbourhoods of $\mu$ such that
 $$\|u(t+\theta_{2})-u(\mu)\|\leq \frac{\varepsilon}{2}, \|u(t+\theta_{1})-u(\mu)\|\leq \frac{\varepsilon}{2}.$$
So, in this case the required equicontinuity follows from the continuity of  $u(t)$ at  $\mu$.

 3. $t+\theta_{1}\geq\mu$. In this case the semigroup  $S(t)$, $t>0$ is compact. So,  $S(t)$ is continuous in the uniform operator topology, and  $S(t)$ is uniformly continuous on $[\mu,T]$. Then
we get
 \begin{equation}\label{37}
\begin{array}{l}
\|u(t+\theta_{1})-u(t+\theta_{2})\|\leq \|S(t+\theta_{1})-S(t+\theta_{1}+r)\|\times
\\
\\\times(\|\varphi(0)\|+\|g(0,\varphi\|+\|g(t+\theta_{1},u_{t+\theta_{1}})-g(t+\theta_{1}+r,u_{t+\theta_{1}+r})\|)+
\\
\\
\int_{0}^{t+\theta_{1}}\|A(S(t+\theta_{1}+r-s)-S(t+\theta_{1}-s))g(s,u_{s})\|ds+\int_{t+\theta_{1}}^{t+\theta_{1}+r}\|A(S(t+\theta_{1}+r-s))g(s,u_{s})\|ds+
\\
\\
+\int_{0}^{t+\theta_{1}}\|(S(t+\theta_{1}+r-s)-S(t+\theta_{1}-s))\|\|f(s,u_{s})\|ds+
\\
\\
+\int_{t+\theta_{1}}^{t+\theta_{1}+r}\|(S(t+\theta_{1}+r-s))\|\|f(s,u_{s})\|ds=I_{1}+I_{2}+I_{3}+I_{4}+I_{5}.
\end{array}
\end{equation}
Uniform continuity implies that   $I_{1}\rightarrow0, r\rightarrow0$ uniformly over $t+\theta_{1}$. Let us estimate  $I_{2}$:
 \begin{equation}\label{38}
\begin{array}{l}
\|g(t+\theta_{1},u_{t+\theta_{1}})-g(t+\theta_{1}+r,u_{t+\theta_{1}+r})\|\leq \|g(t+\theta_{1},u_{t+\theta_{1}})-g(t+\theta_{1}+r,u_{t+\theta_{1}})\|+
\\
\\
+\|g(t+\theta_{1}+r,u_{t+\theta_{1}})-g(t+\theta_{1}+r,u_{t+\theta_{1}+r})\|\leq M\|u_{t+\theta_{1}}-u_{t+\theta_{1}+r}\|_{C}+
\\
\\
+\|g(t+\theta_{1},u_{t+\theta_{1}})-g(t+\theta_{1}+r,u_{t+\theta_{1}})\|.
\end{array}
\end{equation}
The last summand in \eqref{38} tends to zero as $r\rightarrow0$ uniformly over  $t, \theta_{1}$ and $u_{t+\theta_{1}}$ due to condition (H3,3).

Let us estimate  $\|u_{t+\theta_{1}}-u_{t+\theta_{1}+r}\|_{C}$. We have
\begin{equation}\label{39}
 \|u_{t+\theta_{1}}-u_{t+\theta_{1}+r}\|\leq \sup\limits_{\theta\in[-h,0]}\|u(t+\theta_{1}+\theta)-u(t+\theta_{1}+\theta+r)\|.
\end{equation}
Again, we consider three cases.

 1) $t+\theta_{1}+\theta+r<0$. In this case, uniform convergence to zero of \eqref{39} as  $r\rightarrow0$ follows from the uniform continuity of  $\varphi(t), t\in[-h,0]$.

 2) $t+\theta_{1}+\theta<0$, $t+\theta_{1}+\theta+r>0$. Then uniform convergence to zero of \eqref{39} as  $r\rightarrow0$ follows from the uniform continuity of $\varphi(t)$ and $u(t)$ in a neighbourhood of zero point.

 3) $t+\theta_{1}+\theta>0$. Then

\begin{equation}\label{40}
\sup\limits_{\theta\in[-h,0]}\|u(t+\theta_{1}+\theta)-u(t+\theta_{1}+\theta+r)\|\leq \sup\limits_{t\in[\mu,\tau]}\|u(t)-u(t+r)\|.
\end{equation}

Let us estimate $I_{3}$ in \eqref{37}. As the semigroup  $S(t)$ is analytic, then for every  $x\in X$ we have  $S(t)x\in D(A)$ (\cite{Pazy}, Lemma4.2) and $AS(t_{1})S(t)x=S(t_{1})AS(t)x$ for $t_{1}>0$, $t_{2}>0$, $x\in X$.

Thus, we have
$$
\begin{array}{l}
A(S(t+\theta_{1}+r-s)-S(t+\theta_{1}-s))g(s,u_{s})=A(S(\frac{t+\theta_{1}-s}{2}+r)-S(\frac{t+\theta_{1}-s}{2}))
\\
\\
S(\frac{t+\theta_{1}-s}{2})g(s,u_{s})=(S(\frac{t+\theta_{1}-s}{2}+r)-S(\frac{t+\theta_{1}-s}{2}))AS(\frac{t+\theta_{1}-s}{2})g(s,u_{s}).
\end{array}
$$
Therefore,
$$
\begin{array}{l}
I_{3}\leq \int_{0}^{t+\theta_{1}}\|S(\frac{t+\theta_{1}-s}{2}+r)-S(\frac{t+\theta_{1}-s}{2})\|\|A^{1-\alpha}S(\frac{t+\theta_{1}-s}{2})\|\|A^{\alpha}g(s,u_{s})\|ds\leq
\\
\\
\leq \int_{0}^{t+\theta_{1}}\|S(\frac{t+\theta_{1}-s}{2}+r)-S(\frac{t+\theta_{1}-s}{2})\|\|A^{1-\alpha}S(\frac{t+\theta_{1}-s}{2})\|ds\sup\limits_{s\in[0,\tau)}\|g(s,u_{s}\|.
\end{array}
$$

Making a change of variables $\frac{t+\theta_{1}-s}{2}=s_{1}$ and taking into account inequality \eqref{27}, we get
\begin{equation}\label{41}
I_{3}\leq 2\int_{0}^{\tau}\|S(s_{1}+r)-S(s_{1})\|\|A^{1-\alpha}S(s_{1})\|ds_{1}\cdot C.
\end{equation}

 The continuity of $S(t)$ in the uniform operator topology implies that the function under integral in \eqref{41} tends to zero as $r\rightarrow0$. Lemma 2 guarantees that the right-hand part of  \eqref{41} is an integrable function. So, applying  the  dominated convergence theorem, we get
\begin{equation}\label{42}
\int_{0}^{\tau}\|S(s_{1}+r)-S(s_{1})\|\|A^{1-\alpha}S(s_{1})\|ds_{1}\rightarrow 0, \ r\rightarrow 0.
\end{equation}

Uniform over  $t+\theta_{1}$ convergence of $I_{4}$ to zero is a consequence of  \eqref{27} and  \eqref{21}.

To estimate $I_{5}$ we will use Holder's inequality:
\begin{equation}\label{43}
\begin{array}{l}
I_{5}\leq (\int_{0}^{t+\theta_{1}}\|S(t+\theta_{1}+r-s)-S(t+\theta_{1}-s)\|^{q}ds)^{\frac{1}{q}}(\int_{0}^{\tau}m^{p}(s)ds)^{\frac{1}{p}}\leq
\\
\\
\leq(\int_{0}^{\tau}\|S(s+r)-S(s)\|^{q}ds)^{\frac{1}{q}}(\int_{0}^{\tau}m^{p}(s)ds)^{\frac{1}{p}}\rightarrow 0, r\rightarrow 0.
\end{array}
\end{equation}
Using uniform continuity of $S(t)$ and the Lebesgue  dominated convergence theorem, we get
\begin{equation}\label{44}
I_{5}\leq M r^{\frac{1}{q}}(\int_{0}^{\tau}m^{p}(s)ds)^{\frac{1}{p}}\rightarrow 0, r\rightarrow 0.
\end{equation}

Combining  \eqref{37}-\eqref{44} and using inequality  $M_g<1$, we get
$$
\sup\limits_{t+\theta_{1}\in[\mu,\tau)}\|u(t+\theta_{1})-u(t+\theta_{1}+r)\|\rightarrow 0, r\rightarrow\infty,
$$
which means  the uniform continuity over $\theta$ the family of functions
$$\{u(t+\theta), \theta\in[-h,0], t\in[0,\tau)\}.$$
Thus, the set  $P=\{(t,u_{t}): t\in[t_{0},\tau)\}$ belongs to a compact set in $D$, which is a contradiction with Lemma 3.2. Theorem 2 is proved.

\vspace{3 mm}

4. \underline{Applications.}

a) Functional-differential equations of parabolic type.

Let  $Q$ be a bounded domain in  $R^{d}$ with sufficiently smooth boundary $\partial D$. We consider a symmetric  elliptic operator
\begin{equation}\label{45}
A=A(x)=\sum_{i,j=1}^{d}a_{ij}(x)\frac{\partial}{\partial x_{i}\partial x_{j}}=div(a(x),\nabla),
\end{equation}
where $a_{i,j}$ are Holder continuous with the Holder exponent  $\beta\in(0,1)$, bounded,  and for some $C_{0}$
\begin{equation}\label{46}
\sum_{i,j=1}^{d}a_{ij}\eta_{i}\eta_{j}\geq C_{0}\|\eta\|^{2}, \ \eta\in R^{d}.
\end{equation}
We denote $X=L^{2}(Q)=H$, $D(A)=H^{2}(Q)\bigcap H_{0}^{1}(Q)$.

Consider the following initial boundary-value problem
\begin{equation}\label{47}
\begin{array}{l}
\frac{d}{dt}[u(t,x)+b(t,x,\int_{-h}^{0}\|u(t+\theta)\|d\theta)]=div(a(x),\nabla_{x}u(t,x))+f_{1}(t,x,\int_{-h}^{0}\|u(t+\theta)\|d\theta)),
\\
\\
u(t,x)=\varphi(t,x), t\in [t_{0}-h,t_{0}], x\in Q,
\\
\\
u(t,x)=0, x\in \partial Q, t\in[0,T].
\end{array}
\end{equation}
Here $\varphi(t,\cdot)\in C([0,T],L^{2}(Q)), \  t_{0}\in[0,T]$.

Real-valued functions  $b(t,x,y), f_{1}(t,x,y)$ are given for  $t\in [0,T], x\in Q, y\in [0,l], l>0$.

The set $D\subset [0,T)\times C$ is the set  $\{(t,\varphi): t\in [0,T), \varphi\in U\}$, where  $U$ consists of functions  $\varphi\in C$ such that
$\int_{-h}^{0}\|\varphi(\theta,\cdot)\|d\theta \in (0,l)$, and  $\partial U$ consists of those  $\varphi\in C$ for which either $\int_{-h}^{0}\|\varphi(\theta,\cdot)\|d\theta =l$, or $\varphi(\theta,x)=0$ a.e.
 We assume that the initial function  $\varphi(t,x)$ in \eqref{46} also belongs to  $U$. Then  $\partial D=([0,T]\times\partial U)\cup\{T\}\times\overline{U}$.

It is well-known that  $A^{-1}$ is a compact operator  (\cite{Evans}, Sec.6.2),  eigenvalues  $\lambda_{k}$ of  $A$ are real numbers with  $0>\lambda_{1}\geq \lambda_{2}\geq ...$  (\cite{Henry}, Sec.1.4), and the corresponding semigroup  $S(t)$ is compact for $t>0$ and analytic  \cite{Lor} ($A$ is a sectorial self-adjoint operator).

In the sequel, as in \cite{Gris}, we introduce the interpolation space  $D_{A}(\frac{1}{2},2)=H_{0}^{1}$. According to A.17 in \cite{DZ},  $D_{A}(\frac{1}{2},2)$ is isomorphic to $D((-A))^{\frac{1}{2}}$. It means that  $X_{\frac{1}{2}}=H_{0}^{1}$.

Let us consider  assumptions on  $b(t,x,y)$ and $f_{1}(t,x,y)$.

Assume that  $b(t,x,y)$ is a continuous function, and continuity on $t$ is uniform with respect to  $x$ and $y$.

Assume that there exist constants  $L>0, M_{g}>0$ such that
\begin{equation}\label{48}
\begin{array}{l}
|\nabla_{x}b(t,x,y)|\leq L,
\\
\\
|b(t,x,y_{1})-b(t,x.y_{2})|+|\nabla_{x}b(t,x,y_{1})-\nabla_{x}b(t,x,y_{2})|\leq M_{g}|y_{1}-y_{2}|
\end{array}
\end{equation}
for all $t\in[0,T]$, $x\in Q$, $y_{1}, y_{2} \in [0,l]$, and  $2hLM_{g}^{2}meas(Q)<1$.

 Assume that  $f_{1}(t,x,y)$ is measurable with respect to  $t$ and continuous with respect to $x\in Q$, $y\in[0,l]$, and there exists  $m(t)\in L^{p}(0,T) (p>1)$ such that $|f_{1}(t,x,y)|^{2}\leq m(t)(1+|x|^{p}+|y|^{p})$.

After introducing the mappings
\begin{equation}\label{49}
g(t,\varphi)(x):=b(t,x,\int_{-h}^{0}\|\varphi(\theta,\cdot)\|d\theta),
\end{equation}

\begin{equation}\label{50}
f(t,\varphi)(x):=f_{1}(t,x,\int_{-h}^{0}\|\varphi(\theta,\cdot)\|d\theta),
\end{equation}
problem \eqref{47} can be rewritten in the abstract form  \eqref{1}. Let us verify for \eqref{49} and  \eqref{50} conditions of Theorem 1.

From the definition of the uniform metric in $C$ we get continuity of $f(t,\varphi)$ with respect to $\varphi$. Its measurability with respect to  $t$ is a consequence of $t$-measurability of $f_{1}$.

Further,
$$
\|f(t,\varphi)\|^{2}=\int_{Q}f_{1}^{2}(t,x,\int_{-h}^{0}\|\varphi(\theta,\cdot)\|d\theta)dx\leq \int_{Q}m(t)(1+|x|^{p}+h\sup\limits_{\theta\in [-h,0]}\|\varphi(\theta,\cdot)\|^{2})dx.
$$
So, assumption  (H2) holds.

Let us verify assumtions on $g(t,\varphi)$. Its continuity with respect to $t$, uniformly over $\varphi$, is a consequence of continuity of $b$.

Further,
$$
\|g(t,\varphi)\|_{\frac{1}{2}}^{2}=\|g(t,\varphi)\|_{H_{2}^{1}}^{2}=\int_{Q}|g(t,\varphi)(x)|^{2}dx+\int_{Q}|\nabla_{x}g(t,\varphi)(x)|^{2}dx.
$$
But
$$
\int_{Q}|g(t,\varphi)(x)|^{2}dx=\int_{Q}b^{2}(t,x,\int_{-h}^{0}\|\varphi(\theta,\cdot)d\theta\|)dx<\infty
$$
for  $t\in [0,T]$, $x\in Q$, $\varphi\in U$, due to the continuity of $b$.
$$
\int_{Q}|\nabla_{x}g(t,\varphi)(x)|^{2}dx=\int_{Q}|\nabla_{x}b(t,x,\int_{-h}^{0}\|\varphi(\theta,\cdot)d\theta\|)|^{2}dx<L^{2}meas(Q).
$$
Finally, we get
$$
\begin{array}{l}
\|g(t,\varphi_{1})-g(t,\varphi_{2})\|_{\frac{1}{2}}^{2}=\int_{Q}|(g(t,\varphi_{1})(x)-g(t,\varphi_{2})(x))|^{2}dx+
\\
\\
+\int_{Q}|(\nabla_{x}g(t,\varphi_{1})(x)-\nabla_{x}g(t,\varphi_{2})(x))|^{2}dx=\\
\int_Q|b(t,x,\int_{-h}^0\|\varphi_1(\theta,\cdot)\|d\theta)-b(t,x,\int_{-h}^0\|\varphi_2(\theta,\cdot)\|d\theta)|^2+
\\
+\int_{Q}|(\nabla_{x}(b(t,x,\int_{-h}^{0}\|\varphi(\theta_{1},\cdot)d\theta)-b(t,x,\int_{-h}^{0}\|\varphi(\theta_{2},\cdot)d\theta))|^{2}dx\leq
\\
\\
\leq 2M_{g}^{2}\int_{Q}(\int_{-h}^{0}(\|\varphi_{1}(\theta,\cdot)\|-\|\varphi_{2}(\theta,\cdot)\|)d\theta)^{2}\leq 2M_{g}^{2}\int_{Q}(\int_{-h}^{0}(\|\varphi_{1}(\theta,\cdot)-\varphi_{2}(\theta,\cdot)\|)d\theta)^{2}\leq
\\
\\
\leq 2hM_{g}^{2}meas(Q)\sup\limits_{\theta\in[-h,0]}\|\varphi_{1}(\theta,\cdot)-\varphi_{2}(\theta,\cdot)\|^{2}.
\end{array}
$$
So, all assumptions of Theorems 1,2 are fulfilled.

b) parabolic equation with maximum.

Consider the following initial boundary-value problem
\begin{equation}\label{51}
\begin{array}{l}
\frac{d}{dt}(u(t,x)+b(t,x,\max\limits_{s\in I(t)})\|u(s)\|)=Au+f_{1}(t,x,\max\limits_{s\in I(t)}\|u(s)\|),
\\
\\
u(t,x)=\varphi(t,x), t\in [-h,0], x\in Q,
\\
\\
u(t,x)=0, x\in \partial Q, t\in[0,T].
\end{array}
\end{equation}
Here $I(t)=[\beta(t),\alpha(t)]$, $\beta(t), \alpha(t)$ are continuous functions on $[0,T]$ such that  $\beta(t)\leq\alpha(t)\leq t$ and $\min\limits_{t\in[-h,0]}(\beta(t)-t)=-h$.

Assumptions on $A,b,f_1$ are the same as in the previous example.

Problem  \eqref{51} can be rewritten in the abstract form  \eqref{1} if we put
$$
g(t,x)(x)=b(t,x,\max\limits_{\theta\in [\beta(t)-t,\alpha(t)-t]}\|\varphi(\theta,\cdot)\|),
$$
$$
f(t,x)(x)=f_{1}(t,x,\max\limits_{\theta\in [\beta(t)-t,\alpha(t)-t]}\|\varphi(\theta,\cdot)\|).
$$
The set  $U\subset C$ is the set of functions  $\varphi\in C$ such that $\|\varphi(\theta,\cdot)\|\in (0,l)$, $\theta\in [-h,0]$, and $\partial U$ consists of functions from  $C$ such that $\|\varphi(\theta,\cdot)\|\in(0,l)$ and either there exists a point  $\theta\in [-h,0]$ such that  $\|\varphi(\theta,\cdot)\|=l$ or  $\varphi(\theta,x)=0$ for almost all $x\in Q$.

Using the inequality
$$
|\sup\limits_{\theta}\|\varphi_{1}(\theta,\cdot)\|-\sup\limits_{\theta}\|\varphi_{2}(\theta,\cdot)\||\leq \sup\limits_{\theta}\|\varphi_{1}(\theta,\cdot)-\varphi_{2}(\theta,\cdot)\|,
$$
 we can verify all assumptions of Theorems 1,2.

\renewcommand\refname{ References}

\end{document}